\documentclass[12pt]{amsart}
\usepackage{amsmath}
\usepackage{relsize}
\usepackage{amssymb}
\usepackage{amsthm}
\usepackage{mleftright}
\usepackage{amscd}
\usepackage{color}
\usepackage[latin1]{inputenc}
\usepackage{mathpazo}
\usepackage[scaled=.95]{helvet}
\usepackage{courier}
\usepackage[all]{xy}
\usepackage{tikz} 
\usetikzlibrary{babel} 
\usepackage[hypcap]{caption}
\usepackage{hyperref}
\usetikzlibrary{calc} 
\usetikzlibrary{through} 
\usetikzlibrary{patterns} 
\usetikzlibrary{intersections} 
\usetikzlibrary{matrix} 
\usetikzlibrary{cd} 

\usepackage{booktabs} 
\usepackage{array}    
\usepackage{multirow}

\swapnumbers
\headsep=1cm \numberwithin{equation}{section}
\usepackage[colorinlistoftodos]{todonotes}

\newtheorem{theorem}{Theorem}[section]
\newtheorem{lemma}[theorem]{Lemma}
\newtheorem{proposition}[theorem]{Proposition}
\newtheorem{corollary}[theorem]{Corollary}
\newtheorem{fact}[theorem]{Fact}
\newtheorem{facts}[theorem]{Facts}

\theoremstyle{definition}
\newtheorem{definition}[theorem]{Definition}
\newtheorem{remark}[theorem]{Remark}
\newtheorem{remark and definition}{Remark and Definition}[section]

\newtheorem{remark and notation}{Remark and Notation}[section]
\newtheorem{convention}[theorem]{Convention}

\newtheorem{example}[theorem]{Example}

\newtheorem{conjecture}[theorem]{Conjecture}
\newtheorem{question}[theorem]{Question}

\newtheorem*{agra}{Acknowlegment}

\newcommand\depth{\operatorname{depth}}

\newcommand\Ker{\operatorname{\Ker}}

\newcommand{\knds}{\kern-\nulldelimiterspace}

\begin{document}

\title{Hilbert-Kunz Multiplicity of Fiber Product Rings and Nagata Idealizations}



\author{V. H. Jorge-P\'erez}
\address{Universidade de S{\~a}o Paulo -
ICMC
Caixa Postal 668, 13560-970, S{\~a}o Carlos-SP, Brazil}
\email{vhjperez@icmc.usp.br}
\author{Paulo Martins}
\address{Universidade de S{\~a}o Paulo -
ICMC
Caixa Postal 668, 13560-970, S{\~a}o Carlos-SP, Brazil}
\email{paulomartinsmtm@gmail.com}
\author{J. A. Santos-Lima }

\address{Universidade de S{\~a}o Paulo -
ICMC
Caixa Postal 668, 13560-970, S{\~a}o Carlos-SP, Brazil}
\email{seyalbert@gmail.com}


\thanks{Corresponding author: Paulo Martins}

\begin{abstract} The main purpose of this paper is to provide formulas for the Hilbert-Kunz multiplicity of
fiber product rings and Nagata idealizations. We give explicit formulas for the Hilbert-Kunz multiplicity of a fiber product $R \times_T S$, where $R$, $S$, and $T$ are
Noetherian local rings sharing the same characteristic and residue field. We compute the Hilbert-Kunz multiplicity of a general Nagata idealization $R \ltimes M$, where $M$ is a finitely generated $R$-module. Additionally, we provide examples, structural results and establish new bounds for the Hilbert-Kunz multiplicity.
\end{abstract}

\maketitle

\section{Introduction}\label{intro}

Throughout this paper, $(R,\mathfrak{m},k)$ will denote a Noetherian local ring of prime characteristic $p>0$ and dimension $d:=\dim R>0$. It's worth noting that in this context, the residue field $k$ also shares the characteristic $p$. We will denote the length of an $R$-module by $\ell_R(-)$. The well-known Hilbert-Samuel multiplicity is a classical invariant of a local ring that quantifies singularity, where the minimum value of 1 corresponds to a smooth point. In the exploration of singularities, we naturally investigate equimultiple points or the equisingularity of deformations, concepts extensively studied in commutative algebra, algebraic geometry and singularity theory.

In 1983, Monsky \cite{monsky1} introduced a refined version of multiplicity specific to Noethe\-rian local rings of prime characteristic, termed the Hilbert-Kunz multiplicity. It's worth noting that this notion was also introduced by Kunz \cite{kunz1} in 1969, but Monsky provided a more explicit definition. While this invariant remains largely enigmatic, its study has revealed intriguing connections with the classical multiplicity, albeit often exhibiting more intricate behavior. Unlike the Hilbert-Samuel multiplicity, which is always an integer, the Hilbert-Kunz multiplicity may not be. However, this cha\-racteristic enables us to regard it as a more nuanced invariant. Similar to the classical multipli\-city, it is capable of discerning regular rings. 

\begin{definition}
Let $I$ be an $\mathfrak{m}$-primary ideal of $R$  and let $M$ be a finitely generated $R$-module. The \textit{Hilbert-Kunz multiplicity} $e_{HK}(I, M)$ of $M$ with respect to $I$ is defined as follows:
\begin{align}\label{definitionhk}
e_{HK}(I, M):= \lim_{e \to\infty}\frac{\ell_{R}\left(M/I^{[q]}M\right)}{q^{d}},    
\end{align}
where $q=p^{e}$ and $I^{[q]}$ is the ideal generated by $x^q$, $x \in I$. If $M$ is a finitely generated $R$-module $M/I^{[q]}M$ has finite lenght. For simplicity, we can put $e_{HK}(I):=e_{HK}(I, R)$ and
$e_{HK}(R):=e_{HK}(\mathfrak{m})$, if it doesn't cause confusion. Notice that the limit of the right-hand side of (\ref{definitionhk}) always exists and is a real number; see, e.g., Monsky \cite{monsky1}.
\end{definition}

In \cite{kunz2}, Kunz observed that the number $\ell_e(R)=\ell_R(R/\mathfrak{m}^{[q]})/q^d$ is a reasonable measure for the singularity of the local ring $R$ and gives a numerical measure of the 'singular-ness' of $R$. If $R$ is a regular local ring, then $e_{HK}(R)=1$. Otherwise, it is expected that more larger is the value of $e_{HK}(R)-1$, then more severe is the singularity of $R$. The study of fiber products and idealization rings has recently garnered significant attention (see, for instance, \cite{AAM,Buchsbaumness,freitaslima, freitasperez3, freitasperez1, freitasperez2, nan, nasseh2}). Freitas and Perez \cite{freitasperez2} obtained explicit formulas for the Hilbert multiplicity of graded fiber product rings and Strazzanti and Armengou \cite{straz} obtained, as a consequence of their results for some quadratic quotients of the Rees algebra, a formula for the Hilbert-Kunz multiplicity of Nagata idealizations of the form $R \ltimes I$, where $I$ is an $\mathfrak{m}$-primary ideal of $R$. These results motivated us to obtain explicit formulas for the Hilbert-Kunz multiplicity of a fiber product ring \( R \times_T S \) (Theorem \ref{prop:formulas}) and for a general Nagata idealization \( R \ltimes M \) (Corollary \ref{crl:idealization}). These formulas have yielded new results on the structure of fiber product rings, as well as sharper bounds for the Hilbert-Kunz multiplicity of fiber product rings and Nagata idealizations.

Let us briefly describe the contents of the paper. Section \ref{section2} establishes our terminology and  some preparatory facts that will be needed throughout the paper. Motivated by the formulas of Freitas and Perez \cite{freitasperez2} for the Hilbert multiplicity of graded fiber product rings, in Section \ref{sectionfiber}, we obtain explicit formulas for the Hilbert-Kunz multiplicity of a fiber product ring. We provide a simpler proof of a result established by Freitas and Lima \cite{freitaslima} that a non-trivial general fiber product ring $R \times_T S$ is not regular. Additionally, we present a formula for the Hilbert-Kunz multiplicity of the amalgamated duplication $R \bowtie I$ and examples. In Section \ref{sectionidealization}, the main focus is to obtain an explicit formula for the Hilbert-Kunz multiplicity of the general Nagata idealization $R \ltimes M$, where $M$ is a finitely generated $R$-module. Section \ref{sectionapplication} is devoted to studying the behavior of singularities of fiber product rings using a Hilbert-Kunz multiplicity approach. In the final section, we establish bounds for the Hilbert-Kunz multiplicity of fiber product rings and Nagata idealizations. Furthermore, we extend the formula for the Hilbert-Kunz multiplicity of fiber product rings to the more general case of multi-factor fiber products. As an application of our results, we identify new classes of singular rings that satisfy the celebrated Watanabe-Yoshida conjecture.

\section{Setup and Basic Facts}\label{section2}
We begin by recording some fundamental properties of Hilbert-Kunz multiplicity. We then introduce the definitions of fiber products and idealization rings, along with relevant facts that will be used to explore their Hilbert-Kunz multiplicities in the next sections.
\begin{proposition}[Fundamental properties (cf. \cite{Huneke1,kunz1,kunz2,monsky1,watanabe2})]\label{prop2.1} Let $(R,\mathfrak{m},k)$ be a Noetherian local ring of cha\-racteristic $p>0$. Let $I, J$ be $\mathfrak{m}$-primary ideals and let $M$ be a finitely generated $R$-module. Then the following statements hold:
\begin{itemize}
    \item[(1)] If $I \subseteq J$, then $e_{HK}(I) \geq e_{HK}(J)$. 
    \item[(2)]$e_{HK}(I)\geq e_{HK}(R)\geq 1$.
    \item[(3)]\label{null} $e_{HK}(I, M)\geq 0$ and equality holds if and only if $\dim M< \dim R.$    \item[(4)] \textnormal{(Additivity in short exact sequences)} If $0\rightarrow L\rightarrow M\rightarrow N\rightarrow 0$ is a short exact sequence of finitely generated $R$-modules, then \begin{align*}
    e_{HK}(I, M)=e_{HK}(I, L)+e_{HK}(I, N).    
    \end{align*}
    \item[(5)]\label{item} \textnormal{(Associative Formula)}
    \begin{align*}
        e_{HK}(I, M)= \sum_{ \mathfrak{p} \in \operatorname{Assh}(R)}e_{HK}(I, R/\mathfrak{p})\cdot \ell_{R_{\mathfrak{p}}}(M_{\mathfrak{p}}).
    \end{align*}
    where $\operatorname{Assh}(R)$ denotes the set of prime ideals $\mathfrak{p}$ of $R$ with $\dim R/\mathfrak{p}=\dim R$. In particular, if $R$ is a local domain and $M$ is a torsion-free $R$-module of rank $\operatorname{rk}_R M$
, then
$e_{HK}(I, M)= \operatorname{rk}_R M \cdot e_{HK}(I,R)$ .
\item[(6)] If $R$ is a local domain, then $e_{HK}(I,M)=\operatorname{rk}_R M \cdot e_{HK}(I,R)$.
\item[(7)] If $R$ is regular, then $e_{HK}(I)= \ell_R(R/I)$ and $e_{HK}(R)=1$. 
\item[(8)] If $f:(A,\mathfrak{m}) \to (B,\mathfrak{n})$ is a flat local ring homomorphism such that $B/\mathfrak{m}B$ is a field, then $e_{HK}(I)=e_{HK}(IB)$.
\end{itemize}
\end{proposition}
The study of lower bounds for the Hilbert-Kunz multiplicity has been explored by several authors (see e.g., \cite{watanabe2,watanabe1,watanabe3,CDHZ,AE}). Using an elegant argument, Aberbach and Enescu \cite{AE} resolved the fundamental problem of obtaining bounds that are independent of the characteristic in all dimensions.

We say that a local ring $(A,\mathfrak{m}$) is unmixed (resp. quasi-unmixed) if $\dim \widehat{A}/\mathfrak{p}=\dim \widehat{A}$ for every associated (resp. minimal) prime ideal of $\widehat{A}$. Here, $\widehat{A}$ denotes the $\mathfrak{m}$-adic completion of $A$.
\begin{theorem}[Aberbach-Enescu]\label{teo:enescu} Let $(R,\mathfrak{m},k)$ be an unmixed ring of characteristic $p$ and dimension $d\geq 2$. If $R$ is not regular, then 
\begin{align*}
    e_{HK}(R) \geq 1 + \frac{1}{d(d!(d-1)+1)^d}. 
\end{align*}
    
\end{theorem}
\subsection{Fiber product rings} The fiber product structure is studied in various categories. In the case of Noetherian local rings, some fundamental properties can be found in \cite{AAM}.
\begin{definition}\label{definition:fiberproduct}
Let $(R,\mathfrak{m}_R,k)$, $(S,\mathfrak{m}_S,k)$ and $(T,\mathfrak{m}_T,k)$ be
 Noetherian local rings,  and let $R \stackrel{\pi_R}\twoheadrightarrow T \stackrel{\pi_S}\twoheadleftarrow S$ be
  surjective homomorphisms of rings.  \textit{The fiber product} 
$$
R \times_T S=\{(r,s)\in R\times S  \ \mid \ \pi_R(r)=\pi_S(s)  \},
$$
is a Noetherian local ring with maximal
ideal $\mathfrak{m}=\mathfrak{m}_R \times_{\mathfrak{m}_T}{\mathfrak{m}_S}$, residue field $k$, and it is a
subring of the usual direct product $R\times S$ (see \cite[Lemma 1.2]{AAM}). 
Let $\eta_R: R \times_T S\twoheadrightarrow R$
and $\eta_S:R \times_T S\twoheadrightarrow S$ be the natural
projections $(r,s)\mapsto r$ and $(r,s)\mapsto s$, respectively. Then $R \times_TS$ is represented as a pullback diagram:
\begin{align}\label{dia}
\xymatrix{ R \times_T S \ar[d]_{ \eta_R} \ar[r]^{ \; \eta_S} & S \ar[d]^{\pi_S} \\
 R \ar[r]^{\pi_R} & T}
\end{align}
\end{definition}
\begin{convention}\label{conv}
In general, when we use the symbols $R$, $S$ or $T$, we  assume that $R$, $S$ and $T$ are non-zero Noetherian local rings with the maximal ideal $\mathfrak{m}_R$, $\mathfrak{m}_S$ and $\mathfrak{m}_T$, respectively, and with $k$ being the common residue field. The fiber product ring $R \times_TS$ is non-trivial, i.e., $R \neq T \neq S$, and their maximal ideal is denoted by $\mathfrak{m}:= \mathfrak{m}_R \times_{\mathfrak{m}_T} \mathfrak{m}_S$. Note that each (finitely generated) module over $R$, $S$ or $T$ acquires a canonical structure of (finitely generated) $R \times_T S$-module via Diagram (\ref{dia}).
\end{convention}
We will state known facts about fiber product rings to be used in the next sections. 
\begin{fact}[\cite{AAM}]\label{seq}  There is the short exact sequence of $R\times_TS$-modules 
\begin{align*}
0 \rightarrow R \times_T S \rightarrow R \oplus S \rightarrow T \rightarrow 0
\end{align*}
\end{fact}
\begin{fact}[\cite{AAM}, Lemma 1.5]\label{fact:dim} The following (in)equalitie holds 
\begin{align*}
    \dim (R \times_T S)=\max\lbrace \dim R, \dim S \rbrace \geq \min \lbrace \dim R, \dim S \rbrace \geq \dim T
\end{align*}
\end{fact}
\begin{fact}\label{factiso}
Let $k$ be a field and let $I$ and $J$ be ideals of $k[[x_1,\dots,x_n]]$ and $k[[y_1,\dots,y_m]]$, respectively. Then 
\begin{align*}
\frac{k[[x_1,\dots,x_n]]}{I} \times_k \frac{k[[y_1,\dots,y_m]]}{J} \cong \frac{k[[x_1,\dots,x_n,y_1,\dots,y_m]]}{(I+J+(x_iy_j:1 \leq i \leq n, 1 \leq j \leq m))}.
\end{align*}
\end{fact}
\begin{fact}\textnormal{(\cite[Fact 2.9]{nasseh2})}\label{spec}
Let $P = R \times_k S$ be a non-trivial fiber product. Then 
\begin{align*}
    \operatorname{Spec} P = \lbrace \mathfrak{m}_R \oplus \mathfrak{q} \mid \mathfrak{q} \in \operatorname{Spec} S\rbrace \cup \lbrace \mathfrak{p} \oplus \mathfrak{m}_S \mid \mathfrak{p} \in \operatorname{Spec} R \rbrace.
\end{align*}
Following this notation, we have $\operatorname{ht}(\mathfrak{m}_R \oplus \mathfrak{q})= \operatorname{ht}(\mathfrak{q})$ and $P/(\mathfrak{m}_R \oplus \mathfrak{q}) \cong S/\mathfrak{q}$, so we must then have $\dim (P/(\mathfrak{m}_R \oplus \mathfrak{q})) = \dim (S/\mathfrak{q})$ and similarly for the prime $\mathfrak{p} \oplus \mathfrak{m}_S$. Futhermore, if $\mathfrak{q} \neq \mathfrak{m}_S$, then $P_{\mathfrak{m}_R \oplus \mathfrak{q}} \cong S_{\mathfrak{q}}$, and similarly for $\mathfrak{p} \oplus \mathfrak{m}_S$.
\end{fact}
\subsection{Idealization of a module}
Next construction was introduced by Nagata, and its nice properties were studied by Anderson and Winders \cite{idealization}.
\begin{definition}
Let $R$ be a commutative ring and let $M$ be a $R$-module. Then $R \ltimes M: = R \times M$ with coordinate wise addition and multiplication given by $(r_1,m_1)(r_2,m_2):=(r_1r_2,r_1m_2+r_2m_1)$ is a commutative ring with unity (1,0) called the \textit{idealization} of $M$.
\end{definition}
\begin{remark}
Let \( M \) be a finitely generated \( R \)-module over a \( d \)-dimensional Noetherian local ring \( R \) of prime characteristic \( p > 0 \) with maximal ideal \( \mathfrak{m} \) and residue field \( k \). It is well known that the idealization ring \( R \ltimes M \) is also a \( d \)-dimensional Noetherian local ring of characteristic \( p \) with maximal ideal \( \mathfrak{m} \ltimes M \) and residue field \( k \). Furthermore, we can embed \( R \) into \( R \ltimes M \) via \( r \mapsto (r, 0) \), making \( R \ltimes M \) a local \( R \)-algebra with the same residue field as \( R \). This is discussed in detail in \cite{idealization}.
\end{remark}
We will state two already proved theorems about idealization rings to be used in Section \ref{sectionidealization}.
\begin{theorem}[\cite{idealization}, Theorem 3.1]\label{thmidealization} Let $R$ be a commutative ring, $I$ an ideal of $R$, M an $R$-module and $N$ a submodule of $M$. Then $I \ltimes N$ is an ideal of $R \ltimes M$ if and only if $IM \subseteq N$. When $I \ltimes N$ is an ideal, $M/N$ is an $R/I$-module and $(R \ltimes M)/(I \ltimes M) \cong (R/I) \ltimes (M/N)$. In particular $(R \ltimes M)/(0 \ltimes N) \cong R \ltimes (M/N)$ and therefore $(R \ltimes M)/(0 \ltimes M) \cong R$. So the ideals of $R \ltimes M$ containing $0 \ltimes M$ are of the form $J \ltimes M$ for some ideal $J$ of $R$. 
\end{theorem}
\begin{theorem}[\cite{idealization}, Theorem 3.6]\label{primarity} Let $R$ be a commutative ring and $M$ an $R$-module. Let $I$ be an ideal of $R$ and $N$ a submodule of $M$. Then $I \ltimes N$ is primary if and only if either 
\begin{itemize}
    \item[a)] $N=M$ and $I$ is a primary ideal of $R$ or
    \item[b)] $N \subsetneq M$, $IM \subseteq N$, and $I$ and $N$ are $\mathfrak{p}$-primary where $\mathfrak{p}=\sqrt{I}$
\end{itemize}
In either case, $I \ltimes N$ is $\sqrt{I} \ltimes M$-primary.
\end{theorem}

\section{Hilbert-Kunz Multiplicity of Fiber Product Rings}\label{sectionfiber}
In this section, the main focus is to present a formula for the Hilbert-Kunz multipli\-city of the fiber product ring \( R \times_T S \), where \( R \), \( S \), and \( T \) are Noetherian local rings with the same characteristic \( p>0 \) and residue field. In this section, we consider rings $R$, $S$ and $T$ as in Definition \ref{definition:fiberproduct} and adopt the Convention \ref{conv}. To avoid confusion, we note that in this section, \( R \), \( S \), and \( T \) will also be treated as \( R \times_T S \)-modules. Without loss of generality (see Fact \ref{fact:dim}), we may assume that $\dim R \times_T S= \dim R \geq \dim S \geq \dim T$. 
\begin{lemma}\label{prop:fiberproduct}
Let $R\times_TS$ be a fiber product ring of characteristic $p>0$ and let $I$ be an $\mathfrak{m}$-primary ideal of $R\times_TS$. Then the following equality holds
\begin{align*}
e_{HK}(I)= 
\left\{ \begin{array}{llll} e_{HK}(I, R)+e_{HK}(I,S)-e_{HK}(I, T);&\dim R=\dim S=\dim T\\
 e_{HK}(I, R)+e_{HK}(I,S);&  \dim R=\dim S>\dim T\\
e_{HK}(I, R);& \dim R>\dim S \geq \dim T. 
\end{array} \right.
\end{align*}
\begin{proof}
Using the short exact sequence of $R \times_T S$-modules given in Fact \ref{seq} and the additivity of Hilbert-Kunz multiplicity in short exact sequences, we get
\begin{align*}
e_{HK}(I)=e_{HK}(I, R\oplus S)-e_{HK}(I,T).
\end{align*}
Applying the same argument to the exact sequence $0 \rightarrow R \rightarrow R \oplus S \rightarrow S \rightarrow0$
we obtain $e_{HK}(I,R \oplus S)= e_{HK}(I,R)+e_{HK}(I,S)$ and then 
\begin{align*}
    e_{HK}(I)= e_{HK}(I, R)+e_{HK}(I,S)-e_{HK}(I, T).
\end{align*}
Note that the dimensions of $R$, $S$ and $T$ as $R \times_T S$-modules coincide with their respective Krull dimensions. Finally, by using the formula for $\dim( R \times_T S)$ given in Fact \ref{fact:dim} and item (3) of Proposition \ref{prop2.1} we get the desired formulas, nothing that the components of the fiber product that have dimension less than $\dim R \times_T S$ have Hilbert-Kunz multiplicity with respect to $I$ equal to zero.
\end{proof}
\end{lemma}

\begin{theorem}\label{prop:formulas}
Let $R$, $S$ e $T$ be Noetherian local rings of characteristic $p>0$. Then the following equality holds 
\begin{align*}
e_{HK}(R\times_TS)= 
\left\{ \begin{array}{llll} e_{HK}(R)+e_{HK}(S)-e_{HK}(T);&\dim R=\dim S=\dim T\\
 e_{HK}(R)+e_{HK}(S);&  \dim R=\dim S>\dim T\\
e_{HK}( R);& \dim R>\dim S \geq \dim T. 
\end{array} \right.
\end{align*}
\end{theorem}
\begin{proof}
As an immediate consequence of Lemma  \ref{prop:fiberproduct}, taking $I = \mathfrak{m}$, we obtain the following formulas 
\begin{align}\label{formula}
e_{HK}(R\times_TS)= 
\left\{ \begin{array}{llll} e_{HK}(\mathfrak{m}, R)+e_{HK}(\mathfrak{m},S)-e_{HK}(\mathfrak{m}, T);&\dim R=\dim S=\dim T\\
 e_{HK}(\mathfrak{m}, R)+e_{HK}(\mathfrak{m},S);&  \dim R=\dim S>\dim T\\
e_{HK}(\mathfrak{m}, R);& \dim R>\dim S \geq \dim T. 
\end{array} \right.
\end{align}
Now, we will discuss why we can replace $e_{HK}(\mathfrak{m}, R)$, $e_{HK}(\mathfrak{m}, S)$ and $e_{HK}(\mathfrak{m}, T)$ in (\ref{formula}) by $e_{HK}(R)$, $e_{HK}(S)$ and $e_{HK}(T)$, respectively. It is easy to see that $R \times_T S$ also has characteristic $p$. Note that, in general, if $(A,\mathfrak{m}_A,k)$ and $(B,\mathfrak{m}_B,k)$ are $d$-dimensional Noetherian local rings of characteristic $p>0$ and $f:A \to B$ is a surjective ring homomorphism, then $f(\mathfrak{m}_A)=\mathfrak{m}_B$ and so  $\mathfrak{m}_A^{[q]}B=\mathfrak{m}_B^{[q]}B$, with $q=p^e$, for each positive integer $e$. Also, as $\ell_A(N)=\ell_B(N)$, for all $B$-module $N$ of finite lenght, we obtain: 
\begin{align*}
e_{HK}(\mathfrak{m}_A,B) & =  \lim_{e \to\infty}\frac{\ell_{ A}\left(B/\mathfrak{m}^{[q]}_AB\right)}{q^{d}} \\
 & = \lim_{e \to\infty}\frac{\ell_{B}\left(B/\mathfrak{m}_B^{[q]}\right)}{q^{d}} \\ 
  & = e_{HK}(B). 
\end{align*}
Therefore, using the formula for \(\dim R \times_T S\) from Fact \ref{fact:dim} together with the surjective ring homomorphisms in diagram (\ref{dia}), we may substitute \(e_{HK}(\mathfrak{m}, R)\), \(e_{HK}(\mathfrak{m}, S)\), and \(e_{HK}(\mathfrak{m}, T)\) in (\ref{formula}) with \(e_{HK}(R)\), \(e_{HK}(S)\), and \(e_{HK}(T)\), respectively. Here we observe that it suffices to consider only those components of the fiber product \(R \times_T S\) whose dimensions are equal to \(\dim R \times_T S\).
\end{proof}

The study of structural properties of fiber product rings has been addressed by several authors, as \cite{freitasperez1} and \cite{AAM}. In \cite[Proposition 2.6]{nasseh}, they proved that a non-trivial fiber product of the form $R \times_k S$ is not regular. Freitas and Lima \cite[Theorem 4.7(1)]{freitaslima} extended this result to any non-trivial fiber product ring using a Betti numbers approach. We give a simpler proof that any non-trivial fiber product $R \times_T S$ is not regular. We will make use of the singularness of fiber products in the next sections.

\begin{theorem}\label{theorem:singular}
A non-trivial fiber product $R \times_T S$ is not a domain.
\end{theorem}

\begin{proof}
Since \( \pi_R \) and \( \pi_S \) are surjective and not isomorphisms, we have \( \ker(\pi_R) \neq 0 \) and \( \ker(\pi_S) \neq 0 \). Let \( r_0 \in \ker(\pi_R) \setminus \{0\} \) and \( s_0 \in \ker(\pi_S) \setminus \{0\} \). Then,
\[
\pi_R(r_0) = 0 = \pi_S(0) \Rightarrow (r_0, 0) \in R \times_T S,
\]
\[
\pi_R(0) = 0 = \pi_S(s_0) \Rightarrow (0, s_0) \in R \times_T S.
\]
However,
\(
(r_0, 0) \cdot (0, s_0) = (r_0 \cdot 0, 0 \cdot s_0) = (0, 0),
\)
and both \( (r_0, 0) \neq 0 \) and \( (0, s_0) \neq 0 \). Hence, \( R \times_T S \) is not a domain.
\end{proof}

\begin{corollary}\label{singular}
Let $R$, $S$, and $T$ be Noetherian local rings. Then the fiber product \( R \times_T S \) is not regular.
\end{corollary}

\subsection{Amalgamated duplication }
 
Let $R$ be a Noetherian local ring. If we assume $\pi_R=\pi_S:R\to R/I$ (Diagram (\ref{dia})) are the canonical maps and $I$ is a proper $R$-ideal, then $R\times_{R/I}R = R\bowtie I$, where $R\bowtie I$ is the {\it amalgamated duplication of $R$ along $I$,} introduced by D'Anna \cite{danna}. Remark that, if $R$ is local with maximal ideal $\mathfrak{m}$, then $R\bowtie I$ is a local ring with maximal ideal denoted by $\mathfrak{m}\bowtie I.$ Some structural properties of $R\bowtie I$ are summarized in \cite[Theorem 1.8 and 1.9 - 1.9.4]{AAM}.

We can use Theorem \ref{prop:formulas} to obtain a formula for the Hilbert-Kunz Multiplicity of the amalgamated duplication along an ideal. 
\begin{corollary}\label{cor35}
Let $R$ be a Noetherian local ring of characteristic $p>0$ and let $I$ be a proper ideal of $R$. Then
\begin{align*}
 e_{HK}(R \bowtie I) =\left\{ \begin{array}{llll} 2e_{HK}(R)-e_{HK}(R/I);&\dim R=\dim R/I\\
 2e_{HK}(R);&  \dim R>\dim R/I.
\end{array} \right.   
\end{align*}
In particular, if $I$ contains a non-zero divisor, then $e_{HK}(R \bowtie I)= 2e_{HK}(R)$.
\begin{proof}
The proof follows directly using the equality $R \times_{R/I} R = R \bowtie I$ and the formulas given in Theorem \ref{prop:formulas}.
\end{proof}
\end{corollary}
\begin{remark}
    The formula given in Corollary \ref{cor35} generalizes a formula obtained in \cite[Theorem 2.3 and Remark 2.4]{straz} for the Hilbert-Kunz multiplicity of an amalgamated duplication of $R$ along an $\mathfrak{m}$-primary ideal $I$. In this case $\dim R/I= 0$ and then $e_{HK}(R \bowtie I)= 2 e_{HK}(R)$.
\end{remark}
Since $R \bowtie I =R \times_{R/I} R$, the next result is a consequence of Theorem \ref{theorem:singular}.
\begin{corollary}
Let $R$ be a Noetherian local ring and let $I$ be a proper ideal of $R$. Then $R\bowtie I$ is not regular.
\end{corollary}

To end this section, we will use the formulas given in Theorem \ref{prop:formulas} to calculate the Hilbert-Kunz multiplicity of some fiber product rings.

\begin{example}
To compute the Hilbert--Kunz multiplicity of the fiber product \( R^{(r_1)} \times_k R^{(r_2)} \), where \( R^{(r_1)} \) and \( R^{(r_2)} \) are Veronese subrings of \( R = k[[x_1, \dots, x_d]] \), we use Theorem~\ref{prop:formulas} together with the explicit formula from \cite[(1.1) pg. 53]{watanabe3}. Recall that a Veronese subring of order \( r \) is defined by
$
R^{(r)} = \bigoplus_{n \geq 0} R_{nr},
$
where \( R_{nr} \) denotes the homogeneous component of degree \( nr \).

The Hilbert--Kunz multiplicity of the fiber product is then given by:
\[
e_{HK}\big( R^{(r_1)} \times_k R^{(r_2)} \big) = \frac{1}{r_1} \binom{d + r_1 - 1}{r_1 - 1} + \frac{1}{r_2} \binom{d + r_2 - 1}{r_2 - 1}.
\]

Geometrically, the fiber product \( R^{(r_1)} \times_k R^{(r_2)} \) can be interpreted as the union of two formal affine cones over \( \mathrm{Spec}(k) \), each defined by a Veronese subring of different order. 
\end{example}

\begin{example}\label{example1}
Let $k$ be a field of characteristic $p \geq 3$. Consider the hypersurfaces and  $R=k[[x,y,z]]/(xy+z^{n+1})$ and $S=k[[w,u,v]]/(w^2+uv^2+u^{n-1})$, where $n \in \mathbb{N}$. Then
\begin{align*}
  R \times_k S \cong \frac{k[[x,y,z,w,u,v]]}{(xw,xu,xv,yw,yu,yv,zw,zu,zv,xy+z^{n+1},w^2+uv^2+u^{n-1})},
\end{align*}
by Fact \ref{factiso}. In this case, $\dim R = \dim S > \dim k$ and by using the formulas given in \cite[Example 3.18]{Huneke1}, we have
\begin{align*}
e_{HK}(R \times_k S) & = e_{HK} (R) + e_{HK} (S) = \left( 2- \frac{1}{n+1} \right) + \left( 2 - \frac{1}{4(n-2)} \right) \\
 & = 4 - \frac{5n-7}{4(n+1)(n-2)},
\end{align*}
for all positive integer $n \geq 4$, by Theorem \ref{prop:formulas}.
\end{example}
\begin{example}\label{example}
Let $k$ be a field of characteristic $p>0$. Let $R=k[[x]]$, $S=k[[y]]$ and $T=k[[z]]/(z^n)$, where $n \in \mathbb{N}$ and $n \geq 2$. Then $$R \times_T S \cong \frac{k[[x,y]]}{(x^ny-y^2)},$$ by \cite[Corollary 4.2.7]{celikbas}, taking $A=k[[x,y]]/(x^ny-y^2)$, $\mathfrak{p}=(x^n-y)$ and $\mathfrak{q}=(y)$. In this case, $\dim R = \dim S > \dim T$ and we have
\begin{align*}
e_{HK}\left( \frac{k[[x,y]]}{(x^ny-y^2)} \right) = e_{HK}(k[[x]])+e_{HK}(k[[y]])=2, 
\end{align*}
for all positive integer $n \geq 2$, by Theorem \ref{prop:formulas}.
\end{example}

\section{Hilbert-Kunz multiplicity of Nagata idealizations}\label{sectionidealization}
We will start by proving the following lemma, which will play a crucial role in pro\-ving our main theorem of this section. \begin{lemma}\label{lemma:subset}
Let $(R,\mathfrak{m},k)$ be a Noetherian local ring of characteristic $p>0$ and let $M$ be an $R$-module. If $I$ is an ideal of $R$, $N$ is a submodule of $M$ and $J=I\ltimes N$ is an ideal of $R \ltimes M$, or equivalently $IM \subseteq N$, then  $$J^{[p^e]} = I^{[p^e]}(R \ltimes M),$$ 
for all positive integer $e$.
\begin{proof}
First of all, we will fix $q:=p^e$. Note that, for $(r,m) \in R \ltimes M$ we have that $(r,m)^n=(r^n,nr^{n-1}m)$, for all positive integer $n$. So $$(r,m)^q=(r^q,qr^{q-1}m)=(r^q,p^{e}r^{q-1}m)=(r^q,0),$$
since $R$ has characteristic $p$. Then, the generators of the ideal $J^{[q]}$ are of the form $(r^q,0)$, where $r \in I$, and its elements have the form $\sum_{i=1}^n (s_i,m_i)(r_i^q,0)$ where $r_i \in I$, $s_i \in R$ and $m_i \in M$, for $i=1,\dots,n$. Observe that $$\sum_{i=1}^n (s_i,m_i)(r_i^q,0)=\sum_{i=1}^n (r_i^q,0)(s_i,m_i) =\sum_{i=1}^n r_i^q (s_i,m_i)\in I^{[q]}(R \ltimes M).$$ 
Then $J^{[q]} \subseteq I^{[q]}(R \ltimes M)$. On other hand, the generators of the $R$-module $I^{[q]}(R \ltimes M)$ are of the form $r(s,m)$, where $r \in I^{[q]}$ and $(s,m) \in R \ltimes M$. Writing $r=\sum_{i=1}^n y_i x_i^q$, where $y_i \in R$ and $x_i \in I$, for $i=1,\dots,n$, we have
\begin{align*}
  r(s,m) = \sum_{i=1}^n y_ix_i^q(s,m) = \sum_{i=1}^n (y_is,y_im)(x_i^q,0) \in J^{[q]}.
\end{align*}
Then $I^{[q]}(R \ltimes M) \subseteq J^{[q]}$ and the equality holds.
\end{proof}
\end{lemma}
The next theorem provides, as a consequence, a formula for the Hilbert-Kunz multiplicity of a general Nagata idealization $R \ltimes M$, where $M$ is a finitely generated $R$-module.  
\begin{theorem}\label{teo:}
Let $(R,\mathfrak{m},k)$ be a Noetherian local ring of characteristic $p>0$ and let $M$ be a finitely generated $R$-module. If one of the following conditions holds:
\begin{itemize}
    \item[(1)] $I$ is an $\mathfrak{m}$-primary ideal of $R$ and $J=I\ltimes M$.
    \item[(2)] $I$ is an $\mathfrak{m}$-primary ideal of $R$, $N \subsetneq M$ is a submodule of $M$ which is $\mathfrak{m}$-primary and such that $J=I \ltimes N$ is an ideal of $R \ltimes M$ , or equivalently $IM \subseteq N$.
\end{itemize}
The following equalities hold
\begin{align*}
   e_{HK}(J, R\ltimes M) = e_{HK}(I,R)+e_{HK}(I,M)=\lim_{e \to\infty}\frac{\ell_{R}(R\ltimes M/I^{[q]}\ltimes I^{[q]} M)}{q^{d}}.
\end{align*}
In particular, if $\dim M < \dim R$, then $e_{HK}(J, R\ltimes M) = e_{HK}(I,R)$.
\end{theorem}
\begin{proof}
Note that in both cases the ideal $J$ is $\mathfrak{m} \ltimes M$-primary, by Theorem \ref{primarity}. We will fix $q:=p^e$.
Using the canonical short exact sequence of $R$-modules
\begin{align*}
0 \rightarrow R \rightarrow R \ltimes M \rightarrow M \rightarrow0
\end{align*}
we get $e_{HK}(I,R \ltimes M)= e_{HK}(I,R)+e_{HK}(I,M)$. Observe that $\ell_R(P)=\ell_{R \ltimes M}(P)$, for each $R\ltimes M$-module $P$ of finite lenght. Also, using Lemma \ref{lemma:subset} we have the equility $J^{[q]} = I^{[q]}(R \ltimes M)$ and then 
\begin{align*}
e_{HK}(J, R\ltimes M) &  = \lim_{e \to\infty}\frac{\ell_{R\ltimes M}(R\ltimes M/J^{[q]})}{q^{d}}  
\\ & = \lim_{e \to\infty}\frac{\ell_{R}(R\ltimes M/I^{[q]}(R \ltimes M))}{q^{d}} \\
& = e_{HK}(I,R \ltimes M) \\
& = e_{HK}(I,R)+e_{HK}(I,M).
\end{align*}
Also, by Theorem \ref{thmidealization}, we have the isomorphism  
\begin{align*}
    \frac{R \ltimes M}{I^{[q]}\ltimes I^{[q]} M} \cong \frac{R}{I^{[q]}} \ltimes \frac{M}{I^{[q]} M} .
\end{align*}
 Then, considering the exact sequence of $R$-modules of finite length
\begin{align*}
0 \rightarrow \frac{R}{I^{[q]}} \rightarrow \frac{R}{I^{[q]}} \ltimes \frac{M}{I^{[q]} M} \rightarrow  \frac{M}{I^{[q]} M} \rightarrow 0
\end{align*}
And using length additivity in short exact sequences, we obtain
\begin{align*}
\ell_{R} \left( \frac{R \ltimes M}{I^{[q]} \ltimes I^{[q]}M} \right)= \ell_R \left( \frac{R}{I^{[q]}} \right) + \ell_R \left( \frac{M}{I^{[q]}M} \right). 
\end{align*}
Now, taking limits, we get the equality: 
\begin{align*}
\lim_{e \to\infty}\frac{\ell_{R}(R\ltimes M/I^{[q]}\ltimes I^{[q]} M)}{q^{d}} & = \lim_{e \to\infty}\frac{\ell_{R}(R/I^{[q]})}{q^{d}}+ \lim_{e \to\infty}\frac{\ell_{R}(M/I^{[q]} M)}{q^{d}} \\
& = e_{HK}(I,R)+e_{HK}(I,M).
\end{align*}
\end{proof}
\begin{corollary}\label{crl:idealization}
Let $(R,\mathfrak{m},k)$ be a Noetherian local ring of characteristic $p>0$ and let $M$ be a finitely generated $R$-module. Then 
\begin{itemize}
    \item[(1)] $e_{HK}(R \ltimes M) = e_{HK}(R)+e_{HK}(\mathfrak{m},M)$. In particular, if $\dim M < \dim R$, then $e_{HK}(R\ltimes M) = e_{HK}(R)$.
    \item[(2)] If $I$ is a proper ideal of $R$ and $\dim R/I < \dim R$, then $e_{HK}(R\ltimes I) = 2e_{HK}(R)$.
\end{itemize}
\begin{proof}
\begin{itemize}
    \item[(1)] It follows directly from Theorem \ref{teo:}, by taking $J=\mathfrak{m} \ltimes M$. 
    \item[(2)] Using the short exact sequence $0 \rightarrow I \rightarrow R \rightarrow R/I \rightarrow0$,
    we get $e_{HK}(R)=e_{HK}(\mathfrak{m},I)+e_{HK}(\mathfrak{m},R/I)$. Since $\dim R/I < \dim R$, then $e_{HK}(\mathfrak{m},R/I)=0$ (see Proposition \ref{prop2.1}(3)) and therefore $e_{HK}(R)=e_{HK}(\mathfrak{m},I)$. Utilizing item (1), we obtain the desired equality.  
\end{itemize}
\end{proof}
\end{corollary}
\begin{remark}
The formulas given in Corollary \ref{crl:idealization} generalize a formula obtained in \cite[Theorem 2.3 and Remark 2.4]{straz} for the Hilbert-Kunz multiplicity of an idealization ring $R \ltimes I$, where $I$ is an $\mathfrak{m}$-primary ideal of $R$. In this case $\dim R/I= 0$ and then $e_{HK}(R \ltimes I)= 2 e_{HK}(R)$.
\end{remark}

The next corollary concerns providing formulas for the Hilbert-Kunz multiplicity $e_{HK}(J,R \ltimes M)$, where $J$ is an ideal of $R \ltimes M$ as in Theorem \ref{teo:}, by using structural and homological properties of $M$.
\begin{corollary}\label{betti}
Let $(R,\mathfrak{m},k)$ be a Noetherian local ring of characteristic $p>0$ and let $M$ be a finitely generated $R$-module. Let $I$ be an $\mathfrak{m}$-primary ideal of $R$ and let $J$ be an ideal as in Theorem \ref{teo:}. Then
\begin{itemize}
    \item[(1)] $e_{HK}(I,R) \leq e_{HK}(J, R\ltimes M)\leq (1+\mu(M))e_{HK}(I, R)$.
    \item[(2)] If $\operatorname{pd}_R M =n<\infty$, then $$e_{HK}(J,R \ltimes M) = \left( \sum_{i=0}^n(-1)^i \beta_i^R(M) +1 \right) e_{HK}(I,R),$$ 
where $\beta_i^R(M)$ denotes the i-th Betti number of $M$. 
In particular, if $M$ is a free $R$-module of rank $r$, we should have that $e_{HK}(J,R \ltimes M)=(r+1)e_{HK}(I,R)$.
    \item[(3)] If $R$ is a local domain, then $$e_{HK}(J,R \ltimes M)= (\operatorname{rk}_R (M)+1) e_{HK}(I,R).$$ 
\end{itemize}
\begin{proof}
\begin{itemize}
\item[(1)] Considering a short exact sequence $0 \rightarrow K \rightarrow R^{\mu(M)} \rightarrow M \rightarrow 0$
we get $e_{HK}\left( I,R^{\mu(M)} \right)= \mu(M) e_{HK}(I,R) = e_{HK}(I,K)+e_{HK}(I,M)$. Then $\mu(M)  e_{HK}(I, R)\geq  e_{HK}(I, M)$. From this inequality and by Theorem \ref{teo:}, we obtain the result.
    \item[(2)]  Taking a minimal free-resolution of $M$
\begin{align*}
0 \rightarrow R^{\beta_n^R(M)} \rightarrow \cdots \rightarrow R^{\beta_1^R(M)} \rightarrow R^{\beta_0^R(M)} \rightarrow M \rightarrow 0
\end{align*}
and using the additivity of Hilbert-Kunz multiplicity on short exact sequences of finitely generated $R$-modules, we obtain 
\begin{align*}
e_{HK}(I,M)=\sum_{i=0}^n (-1)^i \beta_i^R(M) e_{HK}(I,R).
\end{align*}
Now, utilizing Theorem \ref{teo:} we obtain the desired equality.
  \item[(3)] The proof follows by Theorem \ref{teo:} and item (6) of Proposition \ref{prop2.1}.
\end{itemize}
\end{proof}
\end{corollary}
\begin{example}\label{example2}
Let $R$ be a Noetherian local ring of dimension $d>0$ and characteristic $p>0$. Then 
\begin{align*}
    e_{HK}\left( \frac{R[x_1,\dots,x_n]}{(x_1,\dots,x_n)^2} \right) = (n+1)e_{HK}(R).
\end{align*}
Indeed, if $M=R^n$
, then $R \ltimes M$ is naturally isomorphic to $R[x_1,\dots,x_n]/(x_1,\dots,x_n)^2$, by \cite[Proposition 2.2]{idealization}. Therefore, the equality follows by the above corollary.
\end{example}

The study of structural properties of idealization rings has been addressed by several authors, as \cite{idealization}, \cite{Buchsbaumness}, and \cite{ideal1}.  The next theorem shows that an idealization ring $R \ltimes M$ is not regular. We will make use of the singularness of idealizations in the next sections.

\begin{theorem}\label{idealizationsing} 
Let $R$ be a Noetherian local ring and $M$ a finitely generated $R$-module. Then $R \ltimes M$ is not regular.
\end{theorem}

\begin{proof}
It suffices to show that $R \ltimes M$ is not a domain. In fact, consider
$x = (0, m)$ and $y = (0, m'), \text{with } m, m' \in M, \ m \neq 0, \ m' \neq 0.$ Then $x \cdot y = (0, m)(0, m') = (0, 0 \cdot m' + 0 \cdot m) = (0, 0).$ Therefore, $R \ltimes M$ has zero divisors and is not a domain.  

\end{proof}

\begin{corollary}
Let $A$ be one of the rings $R \times_T S$, $R \bowtie I$, or $R \ltimes M$, each of which is a three-dimensional unmixed local ring of characteristic $p > 0$. Then $e_{HK}(A) \geq \frac{4}{3}$.
\end{corollary}
\begin{proof}
This follows from Corollary~\ref{singular}, Theorem~\ref{idealizationsing}, and \cite[Theorem 3.1 (1)]{watanabe3}.
\end{proof}

\section{Applications to the structure of singularities}\label{sectionapplication}

In this section, we give some applications of the results obtained in the previous sections, focusing on the study of the singularities of the unmixed fiber product ring \( P = R \times_k S \). We show, for example, that if \( e_{HK}(P) = \frac{8}{3} \), then \( P \) is precisely the gluing of two quadratic hypersurfaces in \( \mathbb{A}^4_k \) along the origin. Moreover, we prove that if
$
e_{HK}(P) \in \left( \frac{8}{3}, \frac{17}{6} \right],
$
then \( P \) is exactly the gluing of two hypersurfaces in \( \mathbb{A}^4_k \) along the origin, at least one of type
$
x^2 + y^2 + z^2 + w^c = 0,
$
for some integer \( c \geq 3 \). These define isolated singularities of type \( A \) (or du Val singularities in dimension two), perturbed by a monomial of order at least \( 3 \).

We also analyze the tangent structure of \( P \) at its closed point. That is, we consider the associated graded ring \( \operatorname{gr}_{\mathfrak{m}}(P) \), which can be interpreted as the first-order approximation of \( P \) near the origin-i.e., its tangent cone. 
\begin{remark}\label{unmixedrem}
     Let $A$ be a Noetherian local ring. Then:
$$A \text{ is unmixed} \iff A \text{ is quasi-unmixed and $\widehat{A}$ satisfies Serre's condition } (S_1).$$ 
\end{remark}
Let $(A,\mathfrak{m},k)$ be any local ring of positive dimension. The \textit{associated graded ring} $\operatorname{gr}_{\mathfrak{m}}(A)$ of $A$ with respect to $\mathfrak{m}$ is defined as follows: 
\begin{align*}
\operatorname{gr}_{\mathfrak{m}}(A) := \bigoplus_{n \geq 0} \mathfrak{m}^n/\mathfrak{m}^{n+1}.
\end{align*}
Then $G=\operatorname{gr}_{\mathfrak{m}}(A)$ is a homogeneous $k$-algebra such that $\mathfrak{M}=G_{+}$ is the unique homogeneous maximal ideal of $G$.
\begin{facts}\label{Grad} \textnormal{(\cite[Remark 2.2]{ACL})}
Consider $P=R \times_k S$ with unique maximal ideal $\mathfrak{m}=\mathfrak{m}_R \oplus \mathfrak{m}_S$. Then
    \begin{itemize}
    \item[(1)] $\widehat{P}= \widehat{R} \times_k \widehat{S}$. 
        \item[(2)] $
\operatorname{gr}_{\mathfrak{m}}(P) \cong \operatorname{gr}_{\mathfrak{m}_R}(R) \times_k \operatorname{gr}_{\mathfrak{m}_S}(S)$ .
    \end{itemize}
\end{facts}
\begin{lemma}\label{lemma:dim}
Let $P = R \times_k S$. If $P$ is quasi-unmixed, then $\dim R = \dim S$.
\end{lemma}
\begin{proof}
Since the Krull dimension is invariant under completion, it is enough to prove that $\dim \widehat{R}=\dim \widehat{S}$. Let $\mathfrak{p}_0 \subseteq \cdots \subseteq \mathfrak{p}_m$ and $\mathfrak{q}_0 \subseteq \cdots \subseteq \mathfrak{q}_n$ be two maximal chains of prime ideals realizing the Krull dimensions of $\widehat{R}$ and $\widehat{S}$, respectively. It is then easy to see that $\mathfrak{p}_0 \in \operatorname{Min}(\widehat{R})$, $\mathfrak{q}_0 \in \operatorname{Min}(\widehat{S})$, $\dim \widehat{R}= \dim \widehat{R}/\mathfrak{p}_0$ and $\dim \widehat{S}= \dim \widehat{S}/\mathfrak{q}_0$. By Fact \ref{spec}, we have: 
$$ \operatorname{Spec} \widehat{P} =\operatorname{Spec} \widehat{R} \times_k \widehat{S} = \lbrace \mathfrak{m}_{\widehat{R}} \oplus \mathfrak{q} \mid \mathfrak{q} \in \operatorname{Spec} \widehat{S}\rbrace \cup \lbrace \mathfrak{p} \oplus \mathfrak{m}_{\widehat{S}} \mid \mathfrak{p} \in \operatorname{Spec} \widehat{R} \rbrace.$$
and then $\mathfrak{m}_{\widehat{R}} \oplus \mathfrak{q}_0$ and $\mathfrak{p}_0 \oplus \mathfrak{m}_{\widehat{S}}$ are in $\operatorname{Min}(P)$. Also, since
$\widehat{P}/ (\mathfrak{p}_0 \oplus \mathfrak{m}_{\widehat{S}}) \cong \widehat{R} / \mathfrak{p}_0$ and $\widehat{P} / (\mathfrak{m}_{\widehat{R}} \oplus \mathfrak{q}_0) \cong \widehat{S} / \mathfrak{q}_0$ and $P$ is quasi-unmixed, we have: 
\begin{align*}
    \dim \widehat{R}= \dim \widehat{R}/\mathfrak{p}_0 =  \dim \widehat{S}/\mathfrak{q}_0 = \dim \widehat{S}.   
    \end{align*}
\end{proof}
\begin{proposition}\label{unmixed}
Let $P = R \times_k S$. Then $P$ is unmixed if and only if $R$ and $S$ are unmixed and $\dim R = \dim S$.
\end{proposition}

\begin{proof} By Remark \ref{unmixedrem}, it suffices to show the following:
\begin{itemize}
    \item[(a)] \( P \) is quasi-unmixed if and only if \( R \) and \( S \) are quasi-unmixed and \( \dim R = \dim S \).
    \item[(b)] Let $\dim R = \dim S$. Then $
P$ satisfies $(S_1)$ if and only if $R$ and $S$ satisfy $(S_1)$.
\end{itemize}

\textit{Proof of (a):} Suppose that \( P \) is quasi-unmixed. By Lemma~\ref{lemma:dim} and Fact \ref{fact:dim}, we have \( \dim P = \dim R = \dim S \). Recall that, we have: 
$$ \operatorname{Spec} \widehat{P} = \lbrace \mathfrak{m}_{\widehat{R}} \oplus \mathfrak{q} \mid \mathfrak{q} \in \operatorname{Spec} \widehat{S}\rbrace \cup \lbrace \mathfrak{p} \oplus \mathfrak{m}_{\widehat{S}} \mid \mathfrak{p} \in \operatorname{Spec} \widehat{R} \rbrace.$$
Thus, the minimal primes of \( \widehat{P} \cong \widehat{R} \times_k \widehat{S} \) are exactly those of the form \( \mathfrak{p} \oplus \mathfrak{m}_{\widehat{S}} \), with \( \mathfrak{p} \in \operatorname{Min}( \widehat{R}) \), and \( \mathfrak{m}_{\widehat{R}} \oplus \mathfrak{q} \), with \( \mathfrak{q} \in \operatorname{Min}(\widehat{S}) \). Since
\begin{align}\label{iso}
\widehat{P}/ (\mathfrak{p} \oplus \mathfrak{m}_{\widehat{S}}) \cong \widehat{R} / \mathfrak{p}, \quad \widehat{P} / (\mathfrak{m}_{\widehat{R}} \oplus \mathfrak{q}) \cong \widehat{S} / \mathfrak{q},
\end{align}
it follows from $P$ being quasi-unmixed the equalities \( \dim( \widehat{R} / \mathfrak{p}) = \dim \widehat{P} = \dim \widehat{R} \) and \( \dim(\widehat{S} / \mathfrak{q}) = \dim \widehat{P} = \dim \widehat{S} \), for all $\mathfrak{p} \in \operatorname{Min} (\widehat{R})$ and $\mathfrak{q} \in \operatorname{Min}(\widehat{S})$. Then, \( R \) and \( S \) are quasi-unmixed.

Conversely, assume that \( R \) and \( S \) are quasi-unmixed with \( \dim R = \dim S = d \). Then for all \( \mathfrak{p} \in \operatorname{Min}(\widehat{R}) \), \( \dim(\widehat{R} / \mathfrak{p}) = d \), and similarly for all \( \mathfrak{q} \in \operatorname{Min}(\widehat{S}) \). Thus, since $\dim \widehat{P}=\dim P = d$ (see Fact \ref{fact:dim}) using the isomorphisms in (\ref{iso}), we see that $\dim (\widehat{P}/\mathfrak{t})=\dim \widehat{P}$ for all $\mathfrak{t} \in \operatorname{Min}(\widehat{P})$, as the minimal primes of $\widehat{P}$ are those of the types $\mathfrak{m}_{\widehat{R}} \oplus \mathfrak{q}$ or $\mathfrak{p} \oplus \mathfrak{m}_{\widehat{S}}$, where $\mathfrak{p} \in \operatorname{Min} (\widehat{R})$ and $\mathfrak{q} \in \operatorname{Min}(\widehat{S})$.

\textit{Proof of (b):} Let $\mathfrak{p} \in \operatorname{Spec} P$, $\mathfrak{p} \neq \mathfrak{m} = \mathfrak{m}_R \oplus \mathfrak{m}_S$. Say $\mathfrak{p}= \mathfrak{m}_R \oplus \mathfrak{p}_S$, where \( \mathfrak{p}_S \in \operatorname{Spec}(S) \) is not maximal. Then $P_{\mathfrak{p}} \cong S_{\mathfrak{p}_S}$, by Fact \ref{spec}. We have
\[
\operatorname{depth}(P_{\mathfrak{p}}) = \operatorname{depth}(S_{\mathfrak{p}_S}) \quad \text{and} \quad \dim(P_{\mathfrak{p}}) = \dim(S_{\mathfrak{p}_S}).
\]
Therefore, in this case, we have  
\[
 P\text{ satisfies } (S_1) \text{ at } \mathfrak{p} \iff\operatorname{depth}(P_{\mathfrak{p}}) \geq \min\{1, \dim P_{\mathfrak{p}}\} \iff   S \text{ satisfies } (S_1) \text{ at } \mathfrak{p}_S.
\]
Similarly, if \( \mathfrak{p} = \mathfrak{p}_R \oplus \mathfrak{m}_S \), where \( \mathfrak{p}_R \in \operatorname{Spec}(R) \) is not maximal. We get
\[
P\text{ satisfies } (S_1) \text{ at } \mathfrak{p} \iff\operatorname{depth}(P_{\mathfrak{p}}) \geq \min\{1, \dim P_{\mathfrak{p}}\} \iff R \text{ satisfies } (S_1) \text{ at } \mathfrak{p}_R.
\]
Now, let $\mathfrak{m}=\mathfrak{m}_R \oplus \mathfrak{m}_S$ be the maximal ideal of $P$. We know that $\dim P = \dim R = \dim S$ (see Fact \ref{fact:dim}) and 
\begin{align*}
    \depth (P) = \min \lbrace \depth R,\depth S,1\rbrace,
\end{align*}
by \cite[Corollary 4.5]{freitaslima}. Thus, it is easy to see that $P$ satisfies $(S_1)$ at $\mathfrak{m}$ if and only if $R$ and $S$ satisfy $(S_1)$ at $\mathfrak{m}_R$ and $\mathfrak{m}_S$, respectively.

\end{proof}

\begin{question}
Let $P = R \times_T S$, where $T \neq k$. Is it true that
\[
P \text{ is unmixed } \iff R \text{ and } S \text{ are unmixed and } \dim R = \dim S?
\]
\end{question}

\begin{theorem}
Let $P = R \times_k S$ be a three-dimensional unmixed local ring of characteristic $p > 0$, where $R$ and $S$ are non-regular rings of characteristic $p$. Then the following statements hold:

\begin{enumerate}
    \item[(1)] $e_{HK}(P) \geq \dfrac{8}{3}$.
    
    \item[(2)] Suppose that $k = \overline{k}$ and $\operatorname{char} k \neq 2$. Then the following conditions are equivalent:
    \begin{enumerate}
        \item[(i)] $e_{HK}(P) = \dfrac{8}{3}$.
        \item[(ii)] $\widehat{P} \cong k[[x_1,\dots,x_4,y_1,\dots,y_4]]/(x_1^2 + x_2^2 + x_3^2 + x_4^2, y_1^2 + y_2^2 + y_3^2 + y_4^2, x_iy_j)$.
    \end{enumerate}
\end{enumerate}
\end{theorem}

\begin{proof}
(1): Since $P$ is a three-dimensional unmixed local ring, by Proposition \ref{unmixed}, both $R$ and $S$ must be three-dimensional and unmixed. Thus, since $R$ and $S$ are non-regular, by \cite[Theorem 3.1(1)]{watanabe3}, we have $e_{HK}(R) \geq \dfrac{4}{3}$ and $e_{HK}(S) \geq \dfrac{4}{3}$. Combining this with Theorem \ref{prop:formulas}, we obtain
$$e_{HK}(P) = e_{HK}(R) + e_{HK}(S) \geq \frac{8}{3}.$$

(2): (i) $\Rightarrow$ (ii) By hypothesis and Theorem \ref{prop:formulas}, since $R$ and $S$ share the same dimension, we have:
$$e_{HK}(P) = e_{HK}(R) + e_{HK}(S) = \frac{8}{3}.$$
From \cite[Theorem 3.1(1)]{watanabe3}, as $R$ and $S$ are non-regular, $e_{HK}(R) \geq \dfrac{4}{3}$ and $e_{HK}(S) \geq \dfrac{4}{3}$. Therefore, the equality $e_{HK}(R) + e_{HK}(S) = \dfrac{8}{3}$ implies that $e_{HK}(R) = \dfrac{4}{3}$ and $e_{HK}(S) = \dfrac{4}{3}$. Now, by \cite[Theorem 3.1(2)]{watanabe3}, this yields
$$\widehat{R} \cong \frac{k[[x_1,\dots,x_4]]}{(x_1^2 + x_2^2 + x_3^2 + x_4^2)} \quad \text{and} \quad \widehat{S} \cong \frac{k[[y_1,\dots,y_4]]}{(y_1^2 + y_2^2 + y_3^2 + y_4^2)}.$$
The desired result then follows from Fact \ref{factiso}.

(i) $\Leftarrow$ (ii) Set $H = k[[x,y,z,w]]/(x^2+y^2+z^2+w^2)$. We must then have that $\widehat{P} \cong H \times_k H$ (see Fact \ref{factiso}). Also, $e_{HK}(H)= \frac{4}{3}$ (see e.g., \cite[Theorem 3.1(2)]{watanabe3}). By the flatness of the canonical map $ f: P \to \widehat{P}$, using Proposition \ref{prop2.1}(8) and Theorem \ref{prop:formulas}, we have 
\begin{align*}
    e_{HK}(P)= e_{HK}(\widehat{P})= 2e_{HK}(H)= \frac{8}{3}.
\end{align*}
\end{proof}

Let $(A,\mathfrak{m},k)$ be a local ring of characteristic $p>0$. We say that \(A\) satisfies condition \((\#)\) if \(\mathfrak{m}\) has a minimal reduction \(J\). This condition will be considered in the next theorem.

\begin{theorem}
Let $P=  R \times_k S$ be a three-dimensional unmixed local ring of characteristic $p>0$. Assume further that $k=\overline{k}$ is algebraically closed, $p \neq 2$ and that both $R$ and $S$ are non-regular of characteristic $p$ satisfying the condition ($\#$). Consider the following assertions:
\begin{enumerate}
    \item $\frac{8  }{3} < e_{HK}(P) \leq \frac{17}{6}$   
    \item Either $\widehat{P}
 \cong k[[x_1,\dots,x_4,y_1,\dots,y_4]]/(x_1^2 + x_2^2 + x_3^2 + x_4^{c_1}, y_1^2 + y_2^2 + y_3^2 + y_4^2, x_iy_j)$ or  $P \cong k[[x_1,\dots,x_4,y_1,\dots,y_4]]/(x_1^2 + x_2^2 + x_3^2 + x_4^{c_1}, y_1^2 + y_2^2 + y_3^2 + y_4^{c_2}, x_iy_j) $  for integers $c_1,c_2 \geq 3$.
 \item Either $\operatorname{gr}_{\mathfrak{m}}(P) \cong k[[x_1,x_2,x_3,y_1,\dots,y_4]]/(x_1^2+x_2^2+x_3^2,y_1^2+y_2^2+y_3^2+y_4^2,x_iy_i)$  or  $\operatorname{gr}_{\mathfrak{m}}(P)  \cong k[x_1,x_2,x_3,y_1,y_2,y_3]/(x_1^2 + x_2^2 + x_3^2, y_1^2 + y_2^2 + y_3^2, x_iy_j)$  
\end{enumerate}
Then (1) $\Rightarrow$ (2) $\Rightarrow$ (3).
Moreover, when (2) holds true, we have
\( e_{HK}(P) \geq \dfrac{17}{6} - \dfrac{2}{3c_1^2} \) or \(    e_{HK}(P) \geq 3 - \dfrac{2}{3c_1^2} - \dfrac{2}{3c_2^2} \), respectively.
\begin{proof}
(1) $\Rightarrow$ (2) Let $\frac{8  }{3} < e_{HK}(P) \leq \frac{17}{6}$. Since $P= R \times_k S$ is a three-dimensional unmixed local ring, by Proposition \ref{unmixed}, both $R$ and $S$ must be three-dimensional and unmixed. By \cite[Theorem 3.1(1)]{watanabe3}, we have $e_{HK}(R) \geq \dfrac{4}{3}$ and $e_{HK}(S) \geq \dfrac{4}{3}$. Combining this with Theorem \ref{prop:formulas} and the inequalities $\frac{8  }{3} < e_{HK}(P) \leq \frac{17}{6}$, without loss of generality, we may assume that 
    $e_{HK}(S)=\frac{4}{3}$ and  $\frac{4}{3} < e_{HK}(R) \leq \frac{3}{2}$ or $\frac{4}{3} < e_{HK}(R),e_{HK}(S) \leq \frac{3}{2}$.
In the first case, we have 
\begin{align}\label{hyper}
R \cong \frac{ k[[x_1,\dots,x_4]]}{(x_1^2 + x_2^2 + x_3^2 + x_4^{c_1})} \quad \text{and} \quad \widehat{S} \cong \frac{k[[y_1,\dots,y_4]]}{(y_1^2 + y_2^2 + y_3^2 + y_4^2)},\end{align}
by \cite[Corollary 3.11]{watanabe3} and \cite[Theorem 3.1(2)]{watanabe3}, respectively. Therefore, the isomorphism for $\widehat{P}$ follows by Fact \ref{factiso}. In the second case, we have:
$$R \cong \frac{ k[[x_1,\dots,x_4]]}{(x_1^2 + x_2^2 + x_3^2 + x_4^{c_1})} \quad \text{and} \quad S \cong \frac{k[[y_1,\dots,y_4]]}{(y_1^2 + y_2^2 + y_3^2 + y_4^{c_2})},$$
by \cite[Corollary 3.11]{watanabe3} and the isomorphism for $P$ follows by Fact \ref{factiso}. 

(2) $\Rightarrow$ (3) Since $\operatorname{gr}_{\mathfrak{m}}(P)  \cong \operatorname{gr}_{\widehat{\mathfrak{m}}}(\widehat{P})$ and the associated graded rings of the hypersurfaces in (\ref{hyper}) are known (see \cite[Theorem 3.1(2) and Corollary 3.11]{watanabe3}) the desired isomorphisms follow from Remark \ref{Grad}(1) and Fact \ref{factiso}.
 
The "moreover" part follows by Theorem \ref{prop:formulas}, combining with Fact \ref{factiso} and the formulas given in \cite[Theorem 3.1]{watanabe3} and \cite[Corollary 3.11]{watanabe3}.
\end{proof} 
\end{theorem}

\begin{proposition}
Let \( R \) and \( S \) be Cohen-Macaulay local rings of characteristic \( p > 0 \), both of dimension 2. Let \( T \) be a one-dimensional Cohen-Macaulay local ring of characteristic $p$, and suppose that the residue field \( k \) is algebraically closed. Then $R \times_T S$ is not of "minimal" Hilbert-Kunz multiplicity, that is, $e_{HK}(R \times_T S)>\frac{e(R \times_T S)+1}{2}$. In other words, the associated graded ring \( \operatorname{gr}_{\mathfrak{m}}(R \times_T S) \) is not isomorphic to the Veronese subring \( k[U,V]^{(e)} \), where \( e = e(R)+e(S)\). 
\end{proposition}

\begin{proof}
First, observe that by \cite[Lemma 2.1]{shirogoto}, the fiber product \( R \times_T S \) is Cohen-Macaulay of dimension 2. Then, $e_{HK}(R \times_T S) \geq \frac{e(R \times_T S)+1}{2}$, by \cite[Theorem 2.5]{watanabe1}. Suppose, for contradiction, that
\[
e_{HK}(R \times_T S) =\frac{e(R \times_T S)+1}{2}
\]
From this equality, by Theorem \ref{prop:formulas} and \cite[Corollary 3.4]{freitasperez2}, we derive
\[
e_{HK}(R) + e_{HK}(S) = \frac{e(R) + e(S) + 1}{2}.
\]
Rewriting this, we get
\(
(2e_{HK}(R) - e(R)) + (2e_{HK}(S) - e(S)) = 1
\). However, this leads to a contradiction as
\(
2e_{HK}(R) - e(R) \geq 1 \) and \( 2e_{HK}(S) - e(S) \geq 1
\), by \cite[Theorem 2.5]{watanabe1} again. Then, we must have that $e_{HK}(R \times_T S) > \frac{e(R \times_T S)+1}{2}$.

The "in other words" part follows directly from the equivalence given in \cite[Theorem 3.1]{watanabe1}. 
\end{proof}

\begin{example} 
Let $X \subset \mathbb{A}^3_k$ and $Y \subset \mathbb{A}^3_k$ be affine algebraic surfaces over an algebraically closed field $k$ of characteristic $p > 0$, and let $Z = X \cap Y$. Suppose that the local rings at a point $x \in Z$,
$$
\mathcal{O}_{X,x} \cong k[X]_x \quad \text{and} \quad \mathcal{O}_{Y,x} \cong k[Y]_x,
$$
are both Cohen-Macaulay of dimension 2, and that
$
T = \mathcal{O}_{Z,x} \cong k[Z]_x
$
is Cohen--Macaulay of dimension 1. Then, the associated graded ring of the local ring of the union $X \cup Y$ at $x$ satisfies:
$$
\operatorname{gr}_{\mathfrak{m}}(\mathcal{O}_{X \cup Y, x}) \not\cong k[U, V]^{(e)} \subseteq k[U, V],
$$
where $e = \operatorname{mult}_x(X) + \operatorname{mult}_x(Y)$, and $\operatorname{mult}_x(-)$ denotes the multiplicity of the corresponding variety at $x$ (i.e., the Hilbert-Samuel multiplicity of $\mathcal{O}_{X,x}$ or $\mathcal{O}_{Y,x}$).
\end{example}

    
\section{Bounds for the Hilbert-Kunz multiplicity of fiber products and idealizations}\label{sectionbounds}

In this section, we establish explicit lower bounds for the Hilbert-Kunz multiplicity of fiber products and idealization rings. These bounds are derived by combining our previous results with key known results from recent literature. Although we focus on specific classes of rings, our lower bounds provide significant improvements over existing bounds. Moreover, our results remain valid without certain restrictive assumptions that were crucial in previous results (see, for instance \cite{AE, CDHZ, manuel}).

Note, for example, that the ring $R \times_T S$ is not necessarily unmixed, even if $R$ and $S$ are unmixed. Indeed, if $R = k[[x, y, z]]/(xy)$ and $S = k[[w]]$, we must then have that $
R \times_k S \cong k[[x, y, z,w]]/(xy, xw, yw,zw),
$ is not quasi-unmixed, as it has components of dimensions 1 and 2 (see Lemma \ref{lemma:dim}). Remarkably, even without the unmixedness assumption required in Theorem \ref{teo:enescu}, our results on fiber products yield substantially sharper lower bounds than those obtained in Theorem \ref{teo:enescu} and \cite[Theorem 3.7]{CDHZ}.

\begin{theorem}\label{thm:fiber-product-HK-bound}
Let $P = R \times_T S$ be a $d$-dimensional fiber product ring, where $R$ and $S$ are both (or at least one is) unmixed local rings of characteristic $p > 0$ with $d := \dim R \geq 2$. Then,
\begin{align*}
e_{HK}(R\times_TS)\Rightarrow  
\left\{ \begin{array}{llll}
 =2;&  \dim R=\dim S>\dim T (R, S\mbox{ regular})\\
 \geq 2 + \frac{1}{d\big(d!(d-1) + 1\big)^d};&  \dim R=\dim S>\dim T (R \mbox{ not-regular})\\
 \geq 2\left(1 + \frac{1}{d\big(d!(d-1) + 1\big)^d}\right);&  \dim R=\dim S>\dim T (R,S\mbox{ non-regular})\\
\geq 1 + \frac{1}{d\big(d!(d-1) + 1\big)^d} ;& \dim R>\dim S \geq \dim T (R \mbox{ non-regular}). 
\end{array} \right.
\end{align*}
\end{theorem}
\begin{proof}
Following by Theorem \ref{prop:formulas} and Theorem \ref{teo:enescu}.
\end{proof}

\begin{theorem}\label{Theo 6.2}
Let $P = R \times_T S$ be a local ring of characteristic $p > 0$ with infinite perfect residue field $k$. Assume that $R$ and $S$ are unmixed local rings of dimension $d \geq 2$ and characteristic $p$, where both $R$ and $S$ are non-regular and $F$-finite. 

Let $(\boldsymbol{x})$ and $(\boldsymbol{y})$ be minimal reductions of $\mathfrak{m}_R$ and $\mathfrak{m}_S$ respectively, and define:
\begin{enumerate}
    \item $\mu_1 := \max\{\mu(I) \mid I \text{ is an ideal of } R/(\boldsymbol{x})\}$ and $\mu_2 := \max\{\mu(I) \mid I \text{ is an ideal of } S/(\boldsymbol{y})\},$
    \item $t_1 := \max\{n \in \mathbb{N} \mid \overline{\mathfrak{m}_R^n} \not\subseteq (\boldsymbol{x})\}$ and $t_2 := \max\{n \in \mathbb{N} \mid \overline{\mathfrak{m}_S^n} \not\subseteq (\boldsymbol{y})\}$.
\end{enumerate}
If $d>\dim T$, then
\[
e_{HK}(P) \geq 2 + \min\left\{
\frac{1}{d!},\;
\frac{\mu_1}{e(R) - \mu_1} \cdot \frac{1}{(\lceil d/t_1 \rceil)^d}
\right\}+\min\left\{
\frac{1}{d!},\;
\frac{\mu_2}{e(S) - \mu_2} \cdot \frac{1}{(\lceil d/t_2 \rceil)^d}
\right\},
\]
\end{theorem}
\begin{proof}
The proof follows by Theorem \ref{prop:formulas} and \cite[Theorem 3.7]{CDHZ}.
\end{proof}



Note that the fiber product $R \times_T S$ is not always Cohen-Macaulay or Gorenstein, even when $R$ and $S$ are Cohen-Macaulay or Gorenstein. For example, consider the rings $R$ and $S$ in Example \ref{example1}, which are both 2-dimensional Cohen-Macaulay (in fact Gorenstein) hypersurfaces. Nevertheless, their fiber product $R \times_k S$
fails to be Cohen-Macaulay (and consequently not Gorenstein), as $\dim(R \times_k S) = 2$, while we have $\text{depth}(R \times_k S) = 1$. What is particularly noteworthy is that our results yield equally strong - and in fact sharper - lower bounds for both fiber products and idealizations, without requiring the Cohen-Macaulay or Gorenstein assumptions needed in \cite[Theorem 1.2]{CDHZ} and \cite[Theorem 2.7]{manuel}. We will make these improved bounds explicit in the following results. .

\begin{theorem}\label{Gorenstein}
Let $P = R \times_T S$ be a local ring of characteristic $p>0$ with infinite perfect residue field $k$. Let $R$ and $S$ be reduced $F$-finite Gorenstein rings of dimension $d>\dim T$ and characteristic $p$, and let $s_1$ and $s_2$ the $F$-signatures\footnote{For the definition of F-signature, consult \cite{Huneke1}} of $R$ and $S$, respectively. Let $(\boldsymbol{x})$ and $(\boldsymbol{y})$ be minimal reductions of $\mathfrak{m}_R$ and $\mathfrak{m}_S$, respectively. For any ideals $I$ and $J$ of $R$ and $S$ (resp.) such that $(\boldsymbol{x}) \subseteq I$ and $(\boldsymbol{y}) \subseteq J$, let $\mu_1=\mu (I/(\boldsymbol{x}))$ and $\mu_2=\mu(J/(\boldsymbol{y}))$, then

\[
e_{HK}(P) \geq 2 + \frac{\mu_1}{e(R)-\mu_1}(1-s_1)+\frac{\mu_2}{e(S)-\mu_2}(1-s_2),
\]
\end{theorem}
\begin{proof}
Following by Theorem \ref{prop:formulas} and by \cite[Theorem 1.2]{CDHZ}.
\end{proof}

\begin{theorem}\label{CM}
Let \( P = R \times_T S \) be a local ring of characteristic \( p > 0 \), and set $
\alpha(-) := \max\left\{ \frac{1}{p^d d!}, \frac{1}{p^d e(-)} \right\}.
$
Then,
\begin{itemize}
    \item[(1)] If \( R \) and \( S \) are non-regular Cohen-Macaulay rings of dimension \( d > \dim T \), then
    \[
    e_{HK}(P) > 2 + \alpha(R) + \alpha(S).
    \]

    \item[(2)] If \( R \) is regular and \( S \) is a non-regular Cohen-Macaulay ring with \( \dim R = \dim S = d > \dim T \), then
    \[
    e_{HK}(P) > 2  + \alpha(S).
    \]

    \item[(3)] If \( R \) is a non-regular Cohen-Macaulay ring with \( \dim R = d > \dim S \geq \dim T \), then
    \[
    e_{HK}(P) > 1 + \alpha(R).
    \]
\end{itemize}
\end{theorem}

\begin{proof}
Following by Theorem \ref{prop:formulas} and \cite[Theorem 2.7]{manuel}.
\end{proof}

\begin{remark}\label{boundRemark} Let $A$ be  a Noetherian local ring of characteristic $p > 0$ with dimension $d\geq 2$ and let $M$ be a finitely generated $A$-module. Suppose that $\Lambda=\{\mathfrak{p}\in \operatorname{Assh}(A) \cap \operatorname{Supp}(M) \mid A/\mathfrak{p} \mbox{ is unmixed not regular}\}$. Then, by Associativite Formula (Proposition \ref{prop2.1}(5)) and by Theorem \ref{teo:enescu}, we get

$$e_{HK}(\mathfrak{m},M) \geq \#(\Lambda)\left(1 + \frac{1}{d\big(d!(d-1) + 1\big)^d}\right).$$ 
In particular, if $M=A$ is unmixed and not regular, we get 
$$e_{HK}(A) \geq \#(\Lambda)\left(1 + \frac{1}{d\big(d!(d-1) + 1\big)^d}\right).$$

This last formula improves the formula given by Theorem \ref{teo:enescu} for rings $A$ having at least two distinct primes $\mathfrak{p}\in \operatorname{Assh}(A)$ such that $A/\mathfrak{p}$ is unmixed and not regular. For example, if $P=R\times_kS$ with $\dim R=\dim S$, it has at least two distinct primes $\mathfrak{p} \oplus \mathfrak{m}_S$ and $\mathfrak{m}_R \oplus \mathfrak{q}$ in $\operatorname{Assh}(P)$, and it suffices to require that both $R/\mathfrak{p}$ and $S/\mathfrak{q}$ are unmixed and not regular.
\end{remark}

\begin{proposition}\label{BoundIdealization}
Let \( R\) be an unmixed local ring of characteristic \( p > 0 \) and dimension \( d \geq 2 \). Suppose that \( R \) is not regular, and let \( M \) be a finitely generated \( R \)-module. Then,
\begin{itemize}
\item[(1)] $
e_{HK}(R \ltimes M) \geq (1 + \#\Lambda) \left(1 + \frac{1}{d\left(d!(d - 1) + 1\right)^d}\right).
$
\item[(2)] If $R$ is a domain,  $e_{HK}(R \ltimes M)\geq (\operatorname{rk}_R (M)+1) \left(1 + \frac{1}{d\left(d!(d - 1) + 1\right)^d}\right).$
\end{itemize}
\end{proposition}
\begin{proof}
   (1) Follows from Corollary \ref{crl:idealization}(1), combined with Theorem \ref{teo:enescu} and Remark \ref{boundRemark}, and (2) follows from Corollary \ref{betti}(3) and Theorem \ref{teo:enescu}.
\end{proof}

\subsection{Multi-factor fiber product rings}

We also consider multi-factor fiber products in the local Noetherian setting, which we now define.

\begin{definition}
\label{def:multi_fiber} Let \( R_1, \dots, R_r \) and \( T \) be Noetherian local rings and let 
$
\pi_i : R_i \to T
$
be surjective homomorphisms of local rings. We define the fiber product of \( R_1, \dots, R_r \) over \( T \) as
\[
R_1 \times_T \cdots \times_T R_r = \{ (a_1, \dots, a_r) \in R_1 \oplus \cdots \oplus R_r \mid \pi_i(a_i) = \pi_j(a_j),\ 1 \leq i,j \leq r \}.
\]
\end{definition}

\begin{remark}\label{iterated-fiber}
The multi-factor fiber product construction coincides with iteratively applying the two-factor fiber product construction to the sequence \( R_1, \dots, R_r \) a total of \( r-1 \) times, that is,
\[
R_1 \times_T \cdots \times_T R_r = ((R_1 \times_T R_2) \times_T \cdots ) \times_T R_r.
\]
\end{remark}

From Remark \ref{iterated-fiber}, we get that the multi-factor fiber product inherits almost all intrinsic properties from the two-factor fiber product of rings. In particular, properties such as Noetherianity, locality, and the behavior of the Krull dimension:
\begin{align}\label{multi:dim}
\dim(R_1 \times_T \cdots \times_T R_r) = \max\{ \dim R_i \mid 1 \leq i \leq r \} \geq \min \{ \dim R_i : 1 \leq i \leq r \} \geq \dim T,
\end{align}
as well as the property that \( R_1 \times_T \cdots \times_T R_r \) is a non-regular ring, among others, are preserved. 

From the (in)equalitie (\ref{multi:dim}), it is easy to see that: If $d=\dim(R_1 \times_T \cdots \times_T R_r)= \dim T$, then $\dim R_i=d$ for $1 \leq i \leq r$. Thus, the following formula for the Hilbert-Kunz multiplicity of the multi-factor fiber product $R_1 \times_T \cdots \times_T R_r$ is obtained by applying Theorem~\ref{prop:formulas} $(r-1)$ times, as follows:

\begin{theorem}\label{thm:multi_fiber_cases}
Let \( P = R_1 \times_T \cdots \times_T R_r \) be the multi-factor fiber product ring, and let \( d= \dim (P)\).
Let \( R_1, \dots, R_r \) and \( T \) be Noetherian local rings of characteristic \( p>0 \).
Then,
\[
e_{HK}(P) =
\begin{cases}
\sum\limits_{\substack{1\leq i \leq r }} e_{HK}(R_i) 
- (r-1) e_{HK}(T), &  d = \dim(T), \\[10pt]
\sum\limits_{\substack{1\leq i \leq r \\ \dim(R_i) = d}} e_{HK}(R_i), &  d > \dim(T).
\end{cases}
\]
\end{theorem}



Similarly to Theorem~\ref{thm:fiber-product-HK-bound}, we have the following generalization. This follows from Theorem~\ref{thm:multi_fiber_cases} by applying Theorem \ref{teo:enescu} a finite number of times.

\begin{corollary}\label{multi-bounds}
Let \( P = R_1 \times_T \cdots \times_T R_r \) be the multi-factor fiber product, and let \(d=\dim (P)\).
Let \( R_1, \dots, R_r \)  and \( T \) be Noetherian local rings of characteristic \( p>0 \) such that the rings \( R_i \) of dimension \( d \) are unmixed and non-regular.
Then,
\[
e_{HK}(P) \geq 
\begin{cases}
r\left(1 + \frac{1}{d\big(d!(d-1) + 1\big)^d}\right) -(r-1)e_{HK}(T); & d = \dim(T), \\[10pt]
t\left(1 + \frac{1}{d\big(d!(d-1) + 1\big)^d}\right); &  d > \dim(T),
\end{cases}
\]
where $t$ denotes the number of indices $i$ for which $d = \dim R_i$, in case that $d > \dim T$.

\end{corollary}

\subsection{Application: Watanabe-Yoshida Conjecture for Fiber Products and Idealization Rings}
In this subsection, we present some immediate results regarding the Watanabe-Yoshida conjecture. Thanks to our results from previous sections, we show that there exists a large class of fiber product rings and Nagata idealizations which the conjecture holds under certain assumptions. We state the conjecture below.

\begin{conjecture}[Watanabe-Yoshida \cite{watanabe3}]\label{conj:watanabe}
Let $d \geq 2$ be an integer and $p > 2$ be a prime number. Set $$A_{p,d} := \frac{\overline{\mathbb{F}_p}[[X_0,X_1,\ldots,X_d]]}{(X_0^2 + \cdots + X_d^2)}.$$

Let $(A,\mathfrak{m},k)$ be a $d$-dimensional unmixed local ring of characteristic $p$ and let $k=\overline{\mathbb{F}_p}$. Then the following statements hold:
\begin{enumerate}
    \item If $A$ is not regular, then $e_{\text{HK}}(A) \geq e_{\text{HK}}(A_{p,d}) \geq 1 + m_d$.
    \item If $e_{\text{HK}}(A) = e_{\text{HK}}(A_{p,d})$, then the $\mathfrak{m}$-adic completion of $A$ is isomorphic to $A_{p,d}$ as local rings.
\end{enumerate}
\end{conjecture}

Here, the numbers $m_d$ are defined by the relation 
\begin{align*}
\sec x + \tan x = 1 + \sum_{d=1}^{\infty} m_d x^d \quad (-\pi/2 < x < \pi/2).
\end{align*}
This conjecture has a positive answer for local rings with dimension $\leq 7$ (see \cite{COX}), but for larger dimensions, the problem remains open. Some other positive answers in certain cases are provided, for instance, in \cite{trivedi,enescu}.

For the remainder of this section, we say that a $d$-dimensional Noetherian local ring $A$ of characteristic $p>0$ satisfies (WY1) provided that $A$ satisfies the strong inequality $e_{HK}(A) \geq e_{HK}(A_{p,d})$ of Watanabe-Yoshida conjecture \ref{conj:watanabe}(1). In our work, we obtained explicit formulas for the Hilbert-Kunz multiplicity of fiber products and idealization rings and proved that these rings are not regular. We use these formulas to obtain new classes of singular rings that satisfy (WY1).

\begin{theorem}\label{fiberneq}
Let \( R \), \( S \), and \( T \) be Noetherian local rings of characteristic \( p > 0 \). Assume that $\dim R = \dim R \times_T S$. If \( R \) satisfies \textnormal{(WY1)}, then the fiber product \( R \times_T S \) satisfies \textnormal{(WY1)}.
\end{theorem}

\begin{proof}
By the explicit formulas given in Theorem \ref{prop:formulas}, it is easy to see that this assertion holds when $\dim R \geq \dim S > \dim T$. Thus, we may assume that $\dim R = \dim S= \dim T$, and then $e_{HK}(R \times_T S)=e_{HK}(R)+e_{HK}(S)- e_{HK}(T)$. It is easy to see that $\pi_S(\mathfrak{m}_S)=\mathfrak{m}_T$ and considering the exact sequence of $S$-modules $$0 \rightarrow \operatorname{Ker} \pi_S \rightarrow S \xrightarrow{\pi_S} T \rightarrow 0,$$ we have:
\begin{align*}
    e_{HK}(S)- e_{HK}(\mathfrak{m}_S,T) & = e_{HK}(\mathfrak{m}_S, \operatorname{Ker} \pi_S) 
\\
e_{HK}(S) - \lim_{e \to \infty} \frac{\ell_S(T/\mathfrak{m}_S^{[q]} T)}{q^{\dim S}} & = e_{HK}(\mathfrak{m}_S, \operatorname{Ker} \pi_S) 
\\
e_{HK}(S)-e_{HK}(T)= e_{HK}(S) - \lim_{e \to \infty} \frac{\ell_T(T/\mathfrak{m}_T^{[q]} T)}{q^{\dim T}} & = e_{HK}(\mathfrak{m}_S, \operatorname{Ker} \pi_S) \geq 0.
\end{align*}
Thus, since $R$ satisfies \textnormal{(WY1)}, we have $e_{HK}(R \times_T S )= e_{HK}(R)+e_{HK}(S)-e_{HK}(T) \geq e_{HK}(A_{p,d})$. 
\end{proof}

\begin{corollary}\label{amalga}
Let $R$ be a Noetherian local ring of characteristic $p>0$ and let $I$ be a proper ideal of $R$. If $R$ satisfies \textnormal{(WY1)}, then the amalgamated duplication $R \bowtie I$ satisfies \textnormal{(WY1)}.
\end{corollary}

\begin{theorem}\label{idealizationneq}
Let $R$ be a Noetherian local ring of characteristic $p>0$ and let $M$ be a finitely generated $R$-module. If $R$ satisfies \textnormal{(WY1)}, then the idealization $R \ltimes M$ satisfies \textnormal{(WY1)}. 
\end{theorem}
\begin{proof}
The assertion follows from the formula given in Corollary \ref{crl:idealization}(1).
\end{proof}

An interesting consequence of Theorems \ref{fiberneq} and \ref{idealizationneq} and Corollary \ref{amalga} is due to a result of Enescu and Shimomoto \cite[Theorem 4.6]{enescu}, which shows that every complete intersection ring (not regular) of dimension $d \geq 2$ and characteristic $p>2$ satisfies \textnormal{(WY1)}. It is important to observe that if $R$ is a complete intersection, the rings $R \times_T S$ and $R \ltimes M$ are not always complete intersections (see Examples \ref{example1} and \ref{example2}). Therefore, the next result provides new classes of non-regular local rings that satisfy the strong inequality of Conjecture \ref{conj:watanabe}.

\begin{corollary}
Let $d \geq 2$ and $p>2$. If $R$ is a $d$-dimensional complete intersection ring of characteristic $p$, not regular, then the following assertions hold:
\begin{enumerate}
    \item If $d=\dim R \times_T S$, where $S,T$ are Noetherian local rings of characteristic $p$, then $e_{HK}(R \times_T S) \geq e_{HK}(A_{p,d})$.
    \item If $I$ is a proper ideal of $R$, then $e_{HK}(R \bowtie I) \geq e_{HK}(A_{p,d})$.
    \item If $M$ is a finitely generated $R$-module, then $e_{HK}(R \ltimes M) \geq e_{HK}(A_{p,d})$. 
\end{enumerate}
\end{corollary}

Another interesting consequence of Theorems \ref{fiberneq} and \ref{idealizationneq} and Corollary \ref{amalga} is due to a result of Aberbach of Enescu \cite[Section 4]{AE1}, which shows that the Watanabe-Yoshida Conjecture \ref{conj:watanabe}(1) holds if the ring $R$ has Hilbert-Samuel multiplicity $e(R) \leq 5$. It is important to observe that even if $R$ satisfies $e(R) \leq 5$, the rings $R \times_T S$ and $R \ltimes M$ may have Hilbert-Samuel multiplicity greater than five in several cases (e.g., if $\dim R = \dim S > \dim T$ and $\dim R = \dim M$, respectively). Therefore, the next result provides new classes of non-regular local rings that satisfy the strong inequality of Conjecture \ref{conj:watanabe}.

\begin{corollary}
Let $d \geq 2$ and $p>2$. If $R$ is a $d$-dimensional Noetherian local ring of characteristic $p$, not regular, with $e(R)\leq 5$, then the following assertions hold:
\begin{enumerate}
    \item If $d=\dim R \times_T S$, where $S,T$ are Noetherian local rings of characteristic $p$, then $e_{HK}(R \times_T S) \geq e_{HK}(A_{p,d})$.
  \item If $I$ is a proper ideal of $R$, then $e_{HK}(R \bowtie I) \geq e_{HK}(A_{p,d})$.
    \item If $M$ is a finitely generated $R$-module, then $e_{HK}(R \ltimes M) \geq e_{HK}(A_{p,d})$. 
\end{enumerate}
\end{corollary}

\begin{agra}
 The first author was supported by Sao Paulo Research Foundation (FAPESP) under grant 2019/21181-0. The second author was supported by Sao Paulo Research Foundation (FAPESP) under grant 2022/12114-0.
\end{agra}

\noindent
\textbf{Competing interests:} The authors declare none.

\end{document}